\documentclass[11pt]{article}
\usepackage[utf8]{inputenc}

\usepackage{amsmath}
\usepackage{amscd}
\usepackage{amssymb}
\usepackage{enumerate}
\usepackage{indentfirst}
\usepackage{latexsym}
\usepackage{multicol}
\usepackage{verbatim}
\usepackage{color}
\usepackage{theorem}
\numberwithin{equation}{section}

\def\sgn{\, \hbox{\rm sgn}\,  }

\def\d{\partial}

\def\eps{\eps}

\newtheorem{lem}{Lemma}

\theoremstyle{plain}
\newtheorem{theorem}{Theorem}[section]

\begin{document}


\centerline{\Large Stagnation, Creation, Breaking}

\medskip

\centerline{\bf Piotr Bogus\l aw Mucha}

\begin{center}

{Institute of Applied Mathematics and Mechanics \\
University of Warsaw}

{Email: p.mucha@mimuw.edu.pl}

\end{center}

\begin{abstract}
This note deals with the mono-dimensional equation: $\d_t u -\d_x L(u_x) =f$ with $L(\cdot)$ merely monotone. The goal is to examine the features of facets -- flat regions of graphs of solutions
appearing as  $L(\cdot)$ suffers jumps. We concentrate on examples explaining strong stability, creation and even breaking of facets.

\end{abstract}

{\it MS Classification:} 35K67, 35R35, 35B99, 74E15


\section{Introduction}

Facets, flat regions of interfaces are  inseparable elements of phenomena of a large class of anisotropic phase 
transition models \cite{bbch,Spohn,cahn-taylor}. They  are natural in the crystal growth theory, but also in other phenomena
where isotropy is broken by external conditions \cite{gurtin,garcke,G-R,meirmanov,taylor0}. Creation of such regions is explained
at the microscopic  level. The anisotropy distinguishes selected directions, and then the energy related
 to the interface is
minimized  over such structure. However  in real processes the structure of  interface may be very rough,
using a common language, it is hard to call it {\it smooth}. Effects of creation and  breaking
are typically present,   in addition one can observe very stable (local) structures. 

Mathematical modeling  of such phenomena is still far  from perfect.
We have some systems  formal derived  from classical theories. Such systems are strongly nonlinear,
 the nonlinearities are so
high they cause nonlocal effects \cite{figalli,fukui,koba-giga,mr-siam}. On other hand the current theory of PDEs is not able to handle all
  mathematical challenges. Thus, it is a perfect area to attempt  new 
mathematical theories. The first step is  to understand key phenomena  from the viewpoint
of mathematical analysis in order to build new sets of analytical techniques. And this note is dedicated for this purpose.

In our considerations we address to  mono-dimensional equations. Obstacles appearing in multidimensional problems are related the most with geometry. It causes
very high complexity. The 1d case  can be viewed as a toy model,
 but we shall remember that it can be adapted for description of two dimensional phenomena \cite{garcke,fa-sta,taylor}. 
Here we  study the following class of  parabolic equations
\begin{equation}\label{i0}
 u_t - \partial_x L(u_x) =f
\end{equation}
on an interval with suitable initial and boundary conditions. Operator $L(\cdot)$ is required to be monotone. 
It is well known
that jumps (points of discontinuity) of $L(\cdot)$ cause unusual structure of graphs of solutions. Namely, it generates  flat regions, 
called by the theory:
{\it facets} -- \cite{ggk,koba-giga}.  The key element is the r.h.s. of (\ref{i0}) -- the external force $f$.  
The mathematical goal of this paper is to look closer at structure of solutions to (\ref{i0}) and  understand interaction between $f$ and  the facets.

The aim of the note is  to illustrate three phenomena: stagnation, creation, 
breaking of flat regions of solutions of (\ref{i0}). In the mathematical 
language our analysis will concern:

\smallskip

\noindent
$\clubsuit$ \  stagnation, strong stability of solutions,  a change of data does not 
change the solution; 

\noindent
$\clubsuit$ \ creations of facets, we  prove that from the mathematical point of view 
such shapes are typical, even locally;

\noindent
$\clubsuit$ \ breaking the facet, this phenomenon is spectacular, it is always
associated with a large external force, it can not be a consequence of initial state.

\smallskip 

The stagnation has been noted in results of \cite{meyer,mucha}. The creation was examined in \cite{tvn,mr-joss}. 
The phenomenon of breaking was a subject of many results. We shall mention here the results from \cite{bnp1,bnp2} concerning the weighted curvature flow in the three dimensional case, but without the 
force.   The case of bending/breaking facets under influence of the external force has been investigated in \cite{ggr,gr1,gr2}.

Mathematical problems stated above are arising not only from modeling of evolution of crystals. Such equations arise also from the image processing
\cite{choksi,g,meyer,moll}. We shall keep in mind that the subject can be viewed as a limit case of many general mathematical investigations of nonlinear parabolic systems \cite{alt-luck,mazon2010,AD,bcn,gz}, too.
In other words we will look at mathematical results as interesting on their own. Aspects of applications here are secondary.

This note is a successive step for analysis of (\ref{i0})-type systems.  The key  difference
between the problem from \cite{tvn,mr-non,mr-joss,mr-cb} and this note, is presence of the external force. Thus, here we continue our investigations from \cite{mucha}, where the stationary problem has been examined.

\section{The issue}

Let us describe the problem and basic assumptions. We will examine  the following class of equations
\begin{equation}\label{i1}
 u_t - \d_x L(u_x) = f \mbox{ in } (0,1) \times (0,T)
\end{equation}
with the Dirichlet boundary data, and suitable  initial state
\begin{equation}\label{i1a}
 u(0,t)=u(1,t)=0 \mbox{ \ \ for \ \ } t \in [0,T] \mbox{ \ \ and \ \ }
u|_{t=0}=u_0 \mbox{ in } (0,1).
\end{equation}
The initial datum will be suitably chosen, and the external force $f$ will be sufficiently smooth. We shall keep in mind we are interested in
description of typical  of solutions, thus the full generality of data is not considered in this note.
General form of $L(\cdot)$ does not allow to pinpoint the most interesting examples of behavior of solutions. From that reason we distinguish two cases:
\begin{equation}\label{i2}
 L_0(p)=\sgn p \mbox{ \ \ \ and \ \ \ } L_1(p)=p+\sgn p.
\end{equation}
Nevertheless all operators prescribed  by $L(\cdot)$ suffering jumps will be
in the scope of main interest. Due to required simplicity we
concentrate on the case with just one jump like for $L_1(\cdot)$. 

Let us list main results of our investigations:

\smallskip 

$\blacklozenge$ \ First,  we discuss  regularity of  solutions. Based on the results from \cite{mr-mmas} we show the existence of  solutions. In addition
we are allowed to control the $L_\infty$ bound of $u_{t}$ as well as $L(u_x)_x$, provided suitable information about the initial state.
\smallskip 

$\blacklozenge$ \ Section 4 is dedicated to strong stability -- stagnation -- of facets. In this part we identify a class of initial data, which are insensitive  to changing of the external force.
These examples suggest a very interesting hypothesis/question:  how the structure of the external force determines the changes of solutions.

\smallskip 

$\blacklozenge$ \ The next result explains the effect of creation of facets. We prove that, except  initial time, all extrema must be realized on 
nontrivial intervals. 

\smallskip 

$\blacklozenge$ \ The last point of our investigations is analysis of the phenomenon on breaking facets. We show it just for  one concrete example, but
these considerations describe a generic case for arbitrary flat regions.

\smallskip

Throughout the paper we are trying to keep the standard notation. 

\section{Existence, uniqueness, regularity}

We start with the basic existence result. The monotonicity of the operator $\d_xL(\d_x \cdot)$ immediately yields the uniqueness of solutions to (\ref{i1})-(\ref{i1a}).
In \cite{mr-mmas,mr-joss} we find the following result specifying regularity of unique solutions to systems of (\ref{i1})-(\ref{i1a}) type.

\begin{theorem}\label{t:1}
 Let $L'(\cdot) \geq 0$, $u_0\in L_1(I))$, $u_{0,x}\in BV(I)$, $u_0|_{\{0,1\}}=0$, $f\in C([0,T];W^1_\infty(I))$.
Then there exists a unique weak solution solution to (\ref{i1})-(\ref{i1a}) such that
\begin{equation}
 u_x \in L_\infty(0,T;BV(I)),\qquad u_t \in L_\infty(\delta,T;L_2(I)) \mbox{ for } \delta >0.
\end{equation}
If in addition $L'(\cdot) \geq d>0$ for a constant $d$, then
$ u_{xx} \in L_\infty(\delta,T;L_2(I)).$
\end{theorem}

Very important information is regularity of studied solutions. We prove that for merely monotone $L(\cdot)$
the time derivative $u_t$ is point-wise bounded in the meaning of the $L_\infty$ space,  provided that initial value of $u_t$
is defined suitably. 

\begin{lem}\label{l:1} Let $L'(\cdot) \geq 0$. 
 Let $\d_t f \in L_1(0,T;L_\infty(I))$, $u_t|_{t=0} \in L_\infty(I)$, then
\begin{equation}\label{r1}
\sup_{t\in [0,T]} \|u_t\|_{L_\infty(I)} \leq C(\|u_t|_{t=0}\|_{L_\infty(I)}+\|f_t\|_{L_1(0,T;L_\infty(I))}).
\end{equation}
\end{lem}

{\bf Proof.}
We put our attention on  formal estimates,  purely analytical investigations shall be done at the level of the approximation, see \cite{mr-mmas}.
Differentiate the equation with respect to $t$
\begin{equation}\label{r2}
 u_{tt} - \d_x ( L'(u_x)u_{tx}) = f_t.
\end{equation}
Note that the differentiation in time does not change the boundary conditions.
The simplest approach is by the Moser technique. We test the solution by $|u_t|^{p-1}u_t$ getting
 $\frac{d}{dt} \|u_t\|_{L_p(I)} \leq \|f_t\|_{L_\infty(I)}.$
In the r.h.s. should appear a constant depending on the measure of set $I$, but in our case it is just $1$. Eventually, we prove
\begin{equation}\label{r4}
 \sup_{t\in[0,T]} \|u_t(\cdot,t)\|_{L_p(I)} \leq  \|u_0|_{t=0} \|_{L_\infty}+ \int_0^T \|f_t\|_{L_\infty(I)}.
\end{equation}
Passing with $p \to \infty$ we get (\ref{r1}). The Lemma is proved.

\rightline{\rightline{$\Box$}}

Here we shall put a question concerning the initial data for $u_t$, this value can be obtained only from the equation, 
it means we have to clarify the meaning
 of 
$$
\d_xL(u_{0,x}).
$$ 
This difficulty is related to the programme realized by the author and coworkers \cite{tvn,mr-joss}, however still
the answer for general system is highly non-predictable. The theory \cite{tvn} explains this case just for 
the mono-dimensional total variation flow, 
namely for $L_0(\cdot)$. Hence here we  point two cases which describe us reasonable classes when $u_t|_{t=0}$ 
is determined.

The first one is for $u_t|_{t=0}=0$, this case holds as  $u_0$ is a stationary solution to the following equation
\begin{equation}\label{r6}
 \d_x L(u_{0,x})=f(0).
\end{equation}
In other words, we have

\begin{lem}
 If $u_0$ fulfills the equation (\ref{r6}), then if we consider the evolutionary 
system (\ref{i1}) with the initial data $u_0$ and $f$ such that $f(0)$ is the r.h.s. of (\ref{r6}), then
$u_t|_{t=0}=0$.
\end{lem}

The proof follows from the definition  and results of \cite{mucha}. Thanks to the uniqueness we are able to study the case of an approximation of the system
such that $f(s)=f(0)$ for $s\in [0,\epsilon)$. Then we are ensured that $u\equiv u_0$ on the interval $[0,\epsilon]$. Thus, its time derivative is
zero as well. In particular it is zero for $t=0$ for all $\epsilon$, passing with $\epsilon$ to zero, we obtain the thesis of Lemma in the general form.

\rightline{\rightline{$\Box$}}

Any general answer of the question of the meaning of $u_t|_{t=0}$ is related to definition of the initial quantity: $\d_x L(u_x) - f$. From 
the considerations from \cite{mucha} we still are not able to determine it in a useful way. However we are allowed to consider a very special case for $f(0)=0$.
Then we meet results from \cite{tvn}, for the total variational flow in the mono-dimensional case. They say that if the initial data belong
to an admissible class (which is dense in a suitable class of regularity) then $u_t|_{t=0} \in L_\infty$. So the case of $L_0$ is done.
In  general  the situation is more complex, the know result describes the idea of the almost classical solutions just to 
convex solutions \cite{mr-mmas}.  We do not consider this case since we except better explanation of this issue in the forthcoming results.
A need of redefining the meaning of almost classical solutions to requirements of inhomogeneous systems is required.

\section{Stagnation}

The first feature, we would like to investigate, is the stagnation of facets. The observation says  the length of the flat region depends
of the total force acting on this part. In other words if the change of the force won't change the average then we do not expect change
of the length of the facet. Our examples show even that there will be no influence on the rest of the solution. Thus we obtain the stagnation effect
for the whole solution.

To illustrate this phenomenon we consider the case $L_1(\cdot)$. The analysis can be generalized on the following case
\begin{equation}\label{s0}
 L_2(p)=\sgn p +L_r(p),
\end{equation}
where $L_r$ is a regular $C^1$ function. In other words $L_2(p)$ may suffer just one jump for $p=0$.

Let us  first consider the  simplest case.
We prove

\begin{lem}\label{{l:3}} Consider $L(\cdot)=L_2(\cdot)$ defined by (\ref{s0}).
 Let $u_0 \equiv 0$, $f$ be a smooth function. Then 
\begin{equation}\label{s1}
 u \equiv 0 \mbox{ on } I \times (0,T),
\end{equation}
provided
\begin{equation}\label{s2}
 \int_0^x f(x',t)dx' \in [-1+c,1+c] \mbox{ \ \  for \ \  } t \in [0,t) \mbox{ \ \ \ for some constant $c$.}
\end{equation}
\end{lem}

{\bf Proof. } We construct the solution. Put
\begin{equation}\label{s3}
 \sigma(x,t)=-c + \int_0^x f(x',t)dx' \in [-1,1] =L_2(0).
\end{equation}
Then $u\equiv 0$, since it fulfills the weak formulation of the problem, which is
\begin{equation}\label{s4}
 (0,\phi)+(\sigma,\phi_x)=(f,\phi) \mbox{ for } \phi \in C^\infty_0(0,1).
\end{equation}

\rightline{\rightline{$\Box$}}

{\bf Remark.} To clarify (\ref{s3}) we recall the definition of  weak solutions to (\ref{i1}). Let
$u\in L_\infty(0,T;L_2(\Omega))$ and $\sigma \in L_1(0,T;L_1(I))$ such that
$$
\sigma(x,t) \in L(u_x(x,t)) \mbox{ \ a.e.,}
$$
where $L(u_x)$ is just the composition of multivalued operators;
and the following identity holds
\begin{equation}\label{s5}
 (u_t,\phi)+(\sigma,\phi_x)=(f,\phi) \mbox{  \ \ in } \mathcal{D}'([0,T))
\end{equation}
for all $\phi \in C^\infty([0,T);C^\infty_0(I))$. The $u$ is a weak solution to (\ref{i1}). Such solutions are unique.

\smallskip 

Next, we point a more complex case for $L_1$.

\begin{lem}\label{l:4}
Consider $L(\cdot)=L_1(\cdot)$ defined in (\ref{i2}).
Let $f\in C([0,T);L_\infty(I))$, $u_0$ be a steady solution to the problem
\begin{equation}\label{s6}
 -\d_x(u_x+\sgn u_x)=f(0)
\end{equation}
and $f(0)\geq 0$. In addition we require that $f(0)$ is chosen in such a way that there exist
points $0<\xi_-<\xi_+<1$ such that
\begin{equation}\label{s7}
 \int_{\xi_-}^{\xi_+} f(s,0)ds=2.
\end{equation}
Then the solution $u$ to the following problem 
\begin{equation}\label{s8}
 u_t-\d_x(u_x+\sgn u_x)=f(t)
\end{equation}
is static, it is
\begin{equation}\label{s9}
 u(x,t)=u_0(x) \mbox{ \ \  for } (x,t) \in I \times [0,T),
\end{equation}
provided $f(\cdot)$ fulfills:

\smallskip

\noindent $\spadesuit$ \ supp $f(t)-f(0) \subset [\xi_-,\xi_+]$ for $t\in [0,T)$;

\smallskip
\noindent $\spadesuit$ \ $\int_{\xi_-}^{\xi_+} f(x,t) dx =2$;

\smallskip
\noindent $\spadesuit$ \ $|\int_{a}^{b} f(x,t) dx | <2$ for all $a,b$ such that $\xi_-<a<b<\xi_+$.

\end{lem}

{\bf Proof.} 
The proof follows from considerations in \cite{mucha}, we  show that for all $t$ the force $f(t)$ generates
the same steady solution, so (\ref{s8}) must hold. A key idea follows from regularity of the solutions. Theorem 1 from \cite{mucha} implies the H\"older inequality of
of $u_x$, so we split our analysis into a monotone part, and there we analyze the heat equation, and at the area where the solution is flat.
So a natural boundary condition at the interface is $u_x=0$. In the studied case the solution is expected to be time independent, so it is 
enough to check whether it fulfills the weak formulation, and then the uniqueness ends our considerations.
We keep in mind that all solutions here are unique. 

In Section 6 we prove a more complex result.  Lemma \ref{l:4} is viewed as a particular 
case of  considerations  for Lemmas \ref{l:6} and \ref{l:7}.

\rightline{\rightline{$\Box$}}

\section{Creation}

Next, we consider the phenomenon of creation of flat region of solution. A result of our analysis will be the rule that if $L(\cdot)$
has a jump for value $\bar p$, then the set $\{ u_x = \bar p\}$ consists of isolated non-degenerated closed intervals. We  illustrate 
this behavior concentrated on the model case $L_2(\cdot)$. Thus we study only $L$ 
with jump for $\bar p=0$. In particular, the presence of the force can not break this rule.

\begin{lem}\label{l:5}
 Let $L'(\cdot) \geq 0$ and $f$ be a smooth function. Then each local extremum of the solution $u(\cdot,t)$ is realized 
over non-degenerated interval for $t>0$.
\end{lem}

{\bf Proof.} For fixed $t\in [0,T]$ the space derivative $u_x(\cdot, t) \in TV(I)$, hence $u \in W^1_\infty(I)$.
Assume that we consider a maximum at point $x_0 \in I$.
 Then since  $u_x$ is a TV function,  one can find a sequence $x_n^-(t)$  and $x_n^+(t)$ such that
\begin{equation}\label{c1}
 x_n^-(t) \to x_0^- \mbox{ \ and \ } x_n^+(t) \to x_0^+,
\end{equation}
such that
\begin{equation}\label{c2}
 u(x_n^-,t),\; u(x_n^+,t)<u(x_0,t) \mbox{ \ \  and \ \ } u_x(x_n^-,t)>0,\quad u_x(x_n^+,t) < 0. 
\end{equation}
The integration of the equation (\ref{i1}) over the interval $[x_n^-(t),x_n^+(t)]\times \{t\}$  yields
\begin{equation}\label{c3}
\left. -L(u_x)\right|^{(x_n^+(t),t)}_{(x_n^-(t),t)}=\int_{x_n^-(t)}^{x_n^+(t)} (f(y,t) - u_t(y,t)) dy.
\end{equation}
The choice of  sequences implies 
\begin{equation}\label{c4}
 u_x(x_n^-(t),t) >0 \mbox{ \ \ and \ \ } u_x(x_n^+(t),t) <0.
\end{equation}
Hence we obtain at the limit  the following inequality
\begin{equation}\label{c5}
 2 \leq (x_0^+ - x_0^-)^{1/2}(\int_{x_0^-}^{x_0^+} |f|^2+|u_t|^2dx)^{1/2} \leq (x_0^+ - x_0^-)^{1/2}
(\|f\|_{L_\infty}(t) + \|u_t\|_{L_2}(t)).
\end{equation}
To justify (\ref{c5}) we note that since (\ref{c2}) holds then
\begin{equation}
 \lim_{n\to \infty} \inf L(u_x(x_n^+,t)) \geq 1 \mbox{ \ \ and \ \  } \lim_{n\to \infty} \sup L(u_x(x_n^-,t)) \leq -1.
\end{equation}
It follows that $x_0^-< x_0^+$ and the maximum must hold on a nontrivial interval. On the other hand we obtain the inequality
on the lower bound of the length of the facet. From (\ref{c5}) we get
\begin{equation}\label{c6}
 x_0^+(t) - x_0^-(t) \geq \left(\frac{2}{\|f\|_{L_\infty}(t) + \|u_t\|_{L_2}(t)}\right)^2.
\end{equation}
By Theorem \ref{t:1}, the denominator is  finite for $t>0$, so $x_0^-\neq x_0^+$.
The lower bound (\ref{c6}) can be improved provided we have better information about $u_t$. 
In case of Lemma \ref{l:1} we find then the following bound
\begin{equation}\label{c6a}
 x_0^+(t) - x_0^-(t) \geq \left(\frac{2}{\|f\|_{L_\infty}(t) + \|u_t\|_{L_\infty}(t)}\right).
\end{equation}

\rightline{\rightline{$\Box$}}

\section{Breaking} 

In this section we  analyze the phenomenon of breaking of facets. From previous considerations we have shown that
such objects are very stable, but, as we will see, not unbreakable. We  restrict our-self to the case of 
operator $L_1$, since we will work on the explicit formula.

Let us say few words about solving the problem $$(L_1(u_x))_x=f$$ with the Dirichlet boundary conditions. The following scheme is a 
consequence of our study in \cite{mucha}. We distinguish two subsets  $D_0=\{u_x=0\}$ and $D_1=\{u_x \neq 0\}$. 
On each connected component of $D_1$ we solve the equation $u_{xx}=f$ with $u_x=0$ at the ends, or with the Dirichlet data if the interval 
touches $0$ or $1$. On parts for $D_0$ we have  the energy constraint, ie. $\int_{I_1} f =2$, where $I_1$ is an interval from $D_0$.
To justify the boundary condition $u_x=0$ we note that the steady solution belongs to $W^2_\infty$, so $u_x$ must be a continuous function, 
hence on the points joining components from $D_0$ and $D_1$ we obtain the constraint $u_x=0$. Briefly, it explains solvability of the stationary problem.

Our procedure is the following. We prescribe two steady states: initial one being  a construction from a constant force; and finial state 
being a broken facet generated by a specially chosen force.
First let us determine the initial state. 

\begin{lem}\label{l:6}
 The solution to the problem
\begin{equation}\label{b1}
 L(u_{x})_x=4 \mbox{ in } (0,1), \qquad u(0)=u(1)=0
\end{equation}
has the following form: it is a convex function, symmetric with respect of axis $\{x=1/2\}$
such that 
\begin{equation}\label{b2}
 u(s)=const. \mbox{ \ \ for \ \ } s\in [1/4,3/4].
\end{equation}
The minimum is realized at the value $-1/32$.
\end{lem}

{\bf Proof.} By the results from \cite{mucha} $u \in W^2_\infty(I)$ and the minimum is realized on interval $[\xi_-,\xi_+]$
$$
2=\int_{\xi_-}^{\xi^+} L(u_x)_x dx=4(\xi_+-\xi_-).
$$
The symmetry implies (\ref{b2}). 

\rightline{\rightline{$\Box$}}

The second step is to determine the final state given 
by a more complex force. We would like to add a force acting only on the facet, but without change of the total energy. Hence 
 we consider the following modification
\begin{equation}\label{b3}
 f_\alpha(x)=4+\alpha[-2\chi_{[3/8,5/8]}+\chi_{[1/4,3/4]}] \mbox{ \ \ with  } \alpha \geq 0.
\end{equation}
Then we investigate the problem
\begin{equation}\label{b4}
 L(u_{x})_x=f_\alpha \mbox{ in } (0,1), \qquad u(0)=u(1)=0.
\end{equation}
It is expected that for large $\alpha$ the solution will have two regions of convexity and one 
of concavity, hence there will be two local minima and one local maximum. Of course, all
inner extrema must be realized over non-degenerate intervals.
Using the approach from \cite{mucha} we are able to solve the equation. 

\begin{lem}\label{l:7}
 Let $\alpha >12$, then the solution to (\ref{b4}) has the following properties:

\smallskip 

\noindent
(i) \ all local interior minima are realized on intervals $[1/4,c]$ and $[d,3/4]$ with $1/2<c<3/8$ and $5/8<d<3/4$;

\smallskip 

\noindent
(ii) \ all local interior maximum is realized on the interval $[e,f]$ with $3/8<e<1/2<f<5/8$;

\smallskip 

\noindent
(iii) \ the solution is symmetric with respect of axis $\{x=1/2\}$.

\end{lem}

{\bf Proof.} The symmetry of the solution is obvious. Ends of facets are obtain from the
energy constraints. Then we show that indeed they fulfill the equation, the uniqueness will 
follows the thesis of Lemma. 
Assume that (i) and (ii) holds.
The energy constraints on facets yield
\begin{equation}\label{b5}
 \int_{1/4}^cf_\alpha dx =2
\mbox{ so  } 
(c-1/4)(4+\alpha)=2, \mbox{ it is } c=\frac{2}{4+\alpha} + \frac{1}{4}.
\end{equation}
Provided $\alpha \geq 12$, we get $c\leq 3/8$.
For the maximum, by the symmetry,  we have
\begin{equation}\label{b6}
 \int_{e}^{1/2} f_\alpha dx = -1, \mbox{ it is } e = \frac{1}{2} + \frac{1}{4-\alpha}.
\end{equation}
The same  for $\alpha \geq 12$, we have $e\geq 3/8$.
Thus, the ends of facets are determined. Now we are checking if really we are able to construct 
such solution. We define the solution  on the interval $[0,1/2]$. 
Then $u$ shall fulfill
\begin{equation}\label{b7}
 \begin{array}{l}
  u_{xx}=f_\alpha \mbox{ \ \ in } (0,1/4), \qquad \quad u(0)=0, u_x(1/4)=0;\\
  u(s)=const. \qquad \mbox{ \ \ for } s\in [1/4,c];\\
  u_{xx}=f_\alpha \mbox{ \ \ in } (c,e), \qquad \quad u_x(c)=u_x(e)=0,\\
  u(s)=const. \qquad \mbox{ \ \ for } s\in [e,1/2].
 \end{array}
\end{equation}
Observe that to solve $(\ref{b7})_3$ the compatibility condition $\int_c^e fdx=0$ must be fulfilled. From the definition
$$
\int_c^e fdx=(3/8-c)(4+\alpha)+(e-3/8)(4-\alpha)=\frac{4+\alpha}{8}-2+\frac{4-\alpha}{8}+1=0,
$$
thus the condition is satisfied. It follows that the solution given by (\ref{b7}) exits, since
$$
1/2<c <3/8<e<1/2, \mbox{ \ \ provided } \alpha > 12.
$$
Hence the proof is done. Note that for $\alpha \leq 12$ there is no effects of breaking
of the wall $[1/4,3/4]$. 

A simpler exercise it to check that 
\begin{equation}\label{b7a}
 u^\alpha \leq u^{\alpha'} \mbox{ for } \alpha \leq \alpha'.
\end{equation}
Let us look at the shape of solution on the interval $[c,e]$. Our analysis
gives  that $c(\alpha) \geq c(\alpha')$ as well as $e(\alpha) \leq e(\alpha')$.
On the other hand the solution is a parabola over the set $[c,3/8]$
and $[3/8,e]$ with the highest order coefficient equal to $4+\alpha$ and $4-\alpha$, receptively.
Since at the touch point the first derivatives must be the same we conclude (\ref{b7a}).

\rightline{\rightline{$\Box$}}

\noindent
{\bf Breaking of facets.}
Having Lemmas \ref{l:6} and \ref{l:7} we are prepared to prove the main result of this note.

\begin{lem}\label{l:8}
 Let $\alpha >12$ and $u_0$ be a steady solution to problem (\ref{b1}).  We analyze the solution to the system
\begin{equation}\label{b8}
 u_t-L(u_x)_x= - f_\alpha (x,t)  \mbox{ in } (0,1)\times (0,\infty), \qquad u(0,t)=u(1,t)=0 \mbox{ with } u|_{t=0}=u_0,
\end{equation}
where
\begin{equation}\label{b9}
 f_\alpha (x,t)= 4 + \min\{t,\alpha\}[-2\chi_{[3/8,5/8]}+\chi_{[1/4,3/4]}].
\end{equation}
The solution $u$ has the following structure 
$$
u(x,t)=u_0(x) \mbox{ \ \ for } t \leq 12.
$$
$$
u(x,t) \to u_\alpha(x) \mbox{ \ \ for } t\to \infty,
$$
where $u_\alpha$ is a steady solution to (\ref{b4}) with force $f_\alpha$.

\end{lem}

{\bf Proof.}  Observe that taking the difference between (\ref{b8}) and (\ref{b4}) we get
\begin{equation}\label{b10}
 (u-u_\alpha)_t - \d_x[L(u_x)-L(u_{\alpha,x})]=0 \mbox{ for } t >\alpha,
\end{equation}
then testing by  $u-u_\alpha$ we obtain
\begin{equation}
 \frac{1}{2}\frac{d}{dt} \int_I (u-u_\alpha)^2dx + \int_I [L(u_x)-L(u_{\alpha,x})](u-u_\alpha)_x dx=0.
\end{equation}
Since $L'(\cdot) \geq 2$
\begin{equation}\label{b11}
 \frac{d}{dt} \int_I (u-u_\alpha)^2dx + 2\int_I (u-u_\alpha)_x^2 dx\leq 0
\end{equation}
which guarantees that $u \to u_\alpha$ exponentially in time in the $L_2$ norm.

As a consequence of Lemma \ref{l:4} we observe stagnation for solutions for small $t$. Next, the solution, by (\ref{b11})
tends to $u_\alpha$. Regularity and interpolation relations  imply that $u(t)$ is close to $u_\alpha$ in the $C^{1+a}$-space with $a<1$. 
Thus, the wall over $[1/4,3/4]$ must break into three pieces (we apply here Lemma \ref{l:5}).

\rightline{\rightline{$\Box$}}

{\bf Remark} The above lemma point a way of breaking the facet, the force must be large enough
to create a local energy with opposite sign to one defined on the original wall. In our case
it is related to the parameter $\alpha$, the largeness here is just restricted by the condition $\alpha > 12$.

\bigskip 

\noindent
{\small {\bf Acknowledgment.} The author wishes to express his gratitude to Piotr Rybka for fruitful discussion of this note. The work has been partly supported by 
the NCN Grant No. 2011/01/B/ST1/01197.}

\end{document}